\documentclass[12pt]{amsart}
\usepackage{amssymb}
\usepackage{amsmath,soul}
\usepackage{graphicx}
\newtheorem{theorem}{Theorem}[section]

\newtheorem{exa}[theorem]{Example}

\theoremstyle{definition}

\newtheorem{definition}{Definition}[section]

\newcommand\sgn{\mbox{sgn}}
\newcommand\qqquad{\qquad\qquad\qquad}
\newcommand\qqqquad{\qqquad\qqquad}

\begin{document}
\title{\large A Short Note on Stationary Distributions of Unichain Markov Decision Processes}
\author{\bf Ronald ORTNER}\thanks{~This work was supported in 
part by the the Austrian Science Fund FWF (S9104-N04 SP4) and the IST Programme of the European 
Community, under the PASCAL Network of Excellence, IST-2002-506778. This publication only reflects 
the authors' views.}
\email{rortner@unileoben.ac.at}
\address{\parbox{1.4\linewidth}{Department Mathematik und Informationstechnolgie\\ Montanuniversit\"at Leoben\\
Franz-Josef-Strasse 18\\8700 Leoben, Austria}}
\begin{abstract}
Dealing with unichain MDPs, we consider stationary distributions of policies that coincide in all but $n$ states.
In these states each policy chooses one of two possible actions. We show that the  stationary distributions of 
$n+1$ such policies uniquely determine the stationary distributions of all other such policies. An explicit 
formula for calculation is given.
\end{abstract} 
\maketitle

\markboth{\sf R. Ortner}{\sf A Short Note on Stationary Distributions of Unichain Markov Decision Processes}
\baselineskip 15pt

\section{Introduction}
\begin{definition}
A \emph{Markov decision process} (MDP) $\mathcal{M}$ on a (finite) set of \emph{states} $S$ with a 
(finite) set of \emph{actions} $A$ available in each state $\in S$ consists of 
\begin{enumerate}
\renewcommand{\theenumi}{\roman{enumi}}
\item an initial distribution $\mu_0$ that specifies the probability of starting in some state in $S$,
\item the transition probabilities $p_a(i,j)$ that specify 
the probability of reaching state $j$ when choosing action $a$ in state $i$, and 
\end{enumerate}

A (stationary) \emph{policy} on $\mathcal M$ is a mapping $\pi: S \to A$.
\end{definition}

Note that each policy $\pi$ induces a Markov chain on $\mathcal{M}$. We are interested in MDPs, where 
in each of the induced Markov chains any state is reachable from any other state.

\begin{definition}
An MDP $\mathcal M$ is called \emph{unichain}, if for each policy $\pi$ the Markov
chain induced by $\pi$ is ergodic, i.e. if the matrix $P=(p_{\pi(i)}(i,j))_{i,j\in S}$ is irreducible.
\end{definition}

It is a well-known fact (cf.\ e.g.\ \cite{keme}, p.130ff) that for an ergodic Markov chain with 
transition matrix $P$ there exists a unique invariant and strictly positive distribution $\mu$, such that 
independent of the initial distribution $\mu_0$ one has $\mu_n=\mu_0\bar{P}_n\to \mu$, 
where $\bar{P}_n=\frac{1}{n}\sum_{j=1}^nP^j$.\footnote{Actually, for aperiodic Markov chains 
one has even $\mu_0P^n\to\mu$, while the convergence behavior of periodic Markov chains 
can be described more precisely. However, for our purposes the stated fact is sufficient.}

\section{Main Theorem and Proof}
Given $n$ policies $\pi_1,\pi_2,\ldots,\pi_n$ we say that another policy $\pi$ is a \textit{combination} of $\pi_1,\pi_2,\ldots,\pi_n$, 
if for each state $s$ one has $\pi(s)=\pi_i(s)$ for some $i$.

\begin{theorem}\label{thm:distr}
Let $\mathcal M$ be a unichain MDP and $\pi_1$, $\pi_2$,\ldots,$\pi_{n+1}$ pairwise distinct
policies on $\mathcal M$ that coincide on all but $n$ states $s_1$, $s_2$, \ldots, $s_n$. In these
states each policy applies one of two possible actions, i.e.\ we assume that for each $i$ and each $j$ either
$\pi_i(s_j)=0$ or $\pi_i(s_j)=1$. 
Then the stationary distributions of all combinations of $\pi_1$, $\pi_2$,\ldots,$\pi_{n+1}$ are uniquely 
determined by the stationary distributions $\mu_i$ of the policies $\pi_i$.\\
More precisely, if we represent each combined policy $\pi$ by the word $\pi(s_1) \pi(s_2) \ldots \pi(s_n)$,
we may assume without loss of generality (by swapping the names of the actions correspondingly) that the 
policy $\pi$ we want to determine is $11\ldots 1$. 
Let $S_n$ be the set of permutations of
the elements $\{1,\ldots,n\}$. Then setting
\begin{eqnarray*} 
\Gamma_k &:=&   \{ \gamma \in S_{n+1} \,|\,  \gamma(k)=n+1 \mbox{ and } \pi_j(s_{\gamma(j)})=0 \mbox{ for all } j\neq k \}
\end{eqnarray*} 
one has for the stationary distribution $\mu$ of $\pi$
\[  \mu(s) = \frac{\sum_{k=1}^{n+1} \sum_{\gamma\in\Gamma_k} \sgn(\gamma) \, 
			\mu_k(s) \prod_{j=1\atop j\neq k}^{n+1} \mu_j(s_{\gamma(j)})   }
		{\sum_{s'\in S}\sum_{k=1}^{n+1} \sum_{\gamma\in\Gamma_k} \sgn(\gamma) \, 
			\mu_k(s) \prod_{j=1\atop j\neq k}^{n+1} \mu_j(s_{\gamma(j)}) }. \]
\end{theorem}

For clarification of Theorem \ref{thm:distr}, we proceed with an example.
\begin{exa}\rm
Let $\mathcal M$ be a unichain MDP and $\pi_{000}$, $\pi_{010}$, $\pi_{101}$, $\pi_{110}$
policies on $\mathcal M$ whose actions differ only in three states $s_1$, $s_2$ and $s_3$.
The subindices of a policy correspond to the word $\pi(s_1) \pi(s_2) \pi(s_3)$, so that e.g.\
$\pi_{010}(s_1)=\pi_{010}(s_3)=0$ and  $\pi_{010}(s_2)=1$. Now let $\mu_{000}$, $\mu_{010}$, 
$\mu_{101}$, and $\mu_{110}$ be the stationary distributions of the respective policies.
Theorem \ref{thm:distr} tells us that we may calculate the distributions of all other policies
that play in states $s_1$, $s_2$, $s_3$ action 0 or 1 and coincide with the above mentioned 
policies in all other states. In order to calculate e.g.\ the stationary distribution $\mu_{111}$
of policy $\pi_{111}$ in an arbitrary state $s$, we have to calculate the sets $\Gamma_{000}$, $\Gamma_{010}$, 
$\Gamma_{101}$, and $\Gamma_{110}$. This can be done by interpreting the subindices
of our policies as rows of a matrix. In order to obtain $\Gamma_k$ one cancels row $k$
and looks for all possibilities in the remaining matrix to choose three 0s that neither share
a row nor a column:

\begin{tabbing}
oooooooooo\=ooooooooooooooo\=ooooooooooooooo\=ooooooooooooooo\=ooooooooooooooo\=ooooooooooooooo\kill
\>\st{0 0 0}    \>   \ul{0} 0 0  \>  0 \ul{0} 0  \> \ul{0} 0 0  \>  0 0 \ul{0}\\
\>\ul{0} 1 0     \>  \st{0 1 0}   \>  \ul{0} 1 0  \> 0 1 \ul{0} \> \ul{0} 1 0 \\
\>1 \ul{0} 1     \>   1 \ul{0} 1   \>    \st{1 0 1}  \> 1 \ul{0} 1 \> 1 \ul{0} 1\\
\>1 1 \ul{0}  \>   1 1 \ul{0}     \>   1 1 \ul{0}  \>  \st{1 1 0}  \>  \st{1 1 0}
\end{tabbing}

Each of the matrices now corresponds to a permutation in $\Gamma_k$, where $k$ corresponds
to the cancelled row. Thus $\Gamma_{000}$, $\Gamma_{010}$ and $\Gamma_{101}$ contain
only a single permutation, while  $\Gamma_{110}$ contains two. The respective permutation can be
read off each matrix as follows: note for each row one after another the position of the chosen 0,
and choose $n+1$ for the cancelled row. Thus the permutation for the third matrix is $(2,1,4,3)$.
Now for each of the matrices one has a term that consists of four factors (one for each row).
The factor for a row $j$ is $\mu_j(s')$, where $s'=s$ if row $j$ was cancelled (i.e.\ $j=k$), or equals the state
that corresponds to the column of row $j$ in which the 0 was chosen.
Thus for the third matrix above one gets $\mu_{000}(s_2)\mu_{010}(s_1)\mu_{101}(s)\mu_{110}(s_3)$.
Finally, one has to consider the sign for each of the terms which is the sign of the corresponding
permutation. Putting all together, normalizing the output vector and 
abbreviating $a_i:=\mu_{000}(s_i)$, $b_i:=\mu_{010}(s_i)$, $c_i:=\mu_{101}(s_i)$, and  $d_i:=\mu_{110}(s_i)$
one obtains
\[  \mu_{111}(s) = \frac{\mu_{000}(s)\, b_1 c_2 d_3 - a_1 \mu_{010}(s)\, c_2 d_3 - a_2 b_1 \mu_{101}(s)\, d_3 + 
						a_1 b_3 c_2 \mu_{110}(s)- a_3 b_1 c_2 \mu_{110}(s) }
		{ b_1 c_2 d_3 - a_1 c_2 d_3 - a_2 b_1 d_3 + a_1 b_3 c_2 - a_3 b_1 c_2}. \]
\end{exa}

Theorem \ref{thm:distr} can be obtained from the following more general result
where the stationary distribution of a randomized policy is considered.

\begin{theorem}\label{thm:distr2}
Under the assumptions of Theorem \ref{thm:distr}, the stationary distribution $\mu$ of the 
policy $\pi$ that plays in state $s_i$ ($i=1,\ldots,n$) action 0 with probability $\lambda_i\in [0,1]$ and action 1 with probability 
$(1-\lambda_i)$ is given by
\[  \mu(s) = \frac{\sum_{k=1}^{n+1} \sum_{\gamma\in \Gamma'_k} \sgn(\gamma) \, 
				\mu_k(s) \prod_{j=1\atop j\neq k}^{n+1} f(\gamma(j),j)}
		{\sum_{s'\in S}\sum_{k=1}^{n+1} \sum_{\gamma\in \Gamma'_k} \sgn(\gamma) \, 
				\mu_k(s) \prod_{j=1\atop j\neq k}^{n+1} f(\gamma(j),j)}, \]
where $\Gamma'_k:= \{ \gamma \in S_{n+1} \,|\,  \gamma(k)=n+1 \}$ and
\[ f(i,j) := \begin{cases} \lambda_{i}\,\mu_j(i), & \mbox{if }   \pi_j(i)=1    \\ 
			(\lambda_{i}-1)\,\mu_j(i),   & \mbox{if }   \pi_j(i)=0.
 \end{cases}  \]
\end{theorem}

Theorem \ref{thm:distr} follows from Theorem \ref{thm:distr2} by simply setting $\lambda_i=0$ for $i=1,\ldots,n$.

\begin{proof}[Proof of Theorem \ref{thm:distr2}]
Let $S=\{1,2,\ldots,N\}$ and assume that $s_i=i$ for $i=1,2,\ldots,n$. We denote the probabilities associated
with action 0 with $p_{ij}:=p_0(i,j)$ and those of action 1 with $q_{ij}:=p_1(i,j)$. Furthermore, the probabilities
in the states $i=n+1,\ldots,N$, where the policies $\pi_1,\ldots,\pi_{n+1}$ coincide, are written as $p_{ij}:=p_{\pi_k(i)}(i,j)$
as well. Now setting 
\[ \nu_s := \sum_{k=1}^{n+1} \sum_{\gamma\in \Gamma'_k} \sgn(\gamma) \, 
		\mu_k(s) \prod_{j=1\atop j\neq k}^{n+1} f(\gamma(j),j) \]
and $\nu:=(\nu_s)_{s\in S}$ we are going to show that $\nu P_\pi = \nu$,
where $P_\pi$ is the probability matrix of the randomized policy $\pi$. Since the stationary distribution is
unique, normalization of the vector $\nu$ proves the theorem. Now
\begin{eqnarray*}
(\nu P_\pi )_s &=& \sum_{i=1}^n \nu_i \big(\lambda_i p_{is} + (1-\lambda_i) q_{is} \big) + \sum_{i=n+1}^N \nu_i \, p_{is} \\
&=&  \sum_{i=1}^n \sum_{k=1}^{n+1} \sum_{\gamma\in \Gamma'_k} \sgn(\gamma) \, 
				\mu_k(i) \prod_{j=1\atop j\neq k}^{n+1} f(\gamma(j),j) \big(\lambda_i p_{is} + (1-\lambda_i) q_{is} \big) \\
	&& + \sum_{i=n+1}^N\sum_{k=1}^{n+1} \sum_{\gamma\in \Gamma'_k} \sgn(\gamma) \, 
			\mu_k(i) \prod_{j=1\atop j\neq k}^{n+1} f(\gamma(j),j) \, p_{is}. 
\end{eqnarray*}
Since 
\[ \sum_{i=n+1}^N\mu_k(i) \, p_{is} \;=\; \mu_k(s) - \sum_{i: \pi_k(i)=0} \mu_k(i) \, p_{is} -  \sum_{i: \pi_k(i)=1} \mu_k(i) \, q_{is}, \]
this gives
\begin{eqnarray*}
(\nu P_\pi )_s &=& \sum_{k=1}^{n+1} \sum_{\gamma\in \Gamma'_k} \sgn(\gamma) \, 
		\prod_{j=1\atop j\neq k}^{n+1} f(\gamma(j),j) 
 \Big( \sum_{i=1}^n \mu_k(i) \, \big(\lambda_i p_{is} + (1-\lambda_i) q_{is} \big)\\
	&&  \qqqquad\qquad + \mu_k(s) - \sum_{i: \pi_k(i)=0} \mu_k(i) \, p_{is} -  \sum_{i: \pi_k(i)=1} \mu_k(i) \, q_{is} \Big)\\
&=& \nu_s +  \sum_{k=1}^{n+1} \sum_{\gamma\in \Gamma'_k} \sgn(\gamma) \, 
\prod_{j=1\atop j\neq k}^{n+1} f(\gamma(j),j)  \Big( \sum_{i: \pi_k(i)=0} \mu_k(i) \, (\lambda_i-1) (p_{is}  - q_{is}) \\
	&&  \qqqquad\qqquad   +\sum_{i: \pi_k(i)=1} \mu_k(i) \, \lambda_i (p_{is}  - q_{is})  \Big) \\
&=& \nu_s +  \sum_{k=1}^{n+1} \sum_{\gamma\in \Gamma'_k} \sgn(\gamma) \, 
	\prod_{j=1\atop j\neq k}^{n+1} f(\gamma(j),j) \sum_{i=1}^{n}  (p_{is}  - q_{is}) f(i,k)  \\
&=&\nu_s +   \sum_{i=1}^{n}(p_{is}  - q_{is})\sum_{k=1}^{n+1} \sum_{\gamma\in \Gamma'_k}  \,\sgn(\gamma) \, 
f(i,k)\,\prod_{j=1\atop j\neq k}^{n+1} f(\gamma(j),j)
\end{eqnarray*}
Now it is easy to see that  $\sum_{k=1}^{n+1} \sum_{\gamma\in \Gamma'_k}  \,\sgn(\gamma) \, 
f(i,k)\,\prod_{j=1\atop j\neq k}^{n+1} f(\gamma(j),j)=0$: fix $k$ and some permutation $\gamma\in\Gamma_k'$, and 
let $l:=\gamma^{-1}(i)$. 
Then there is exactly one permutation $\gamma'\in\Gamma_l'$, such that $\gamma'(j)=\gamma(j)$ for $j\neq k,l$ and $\gamma'(k)=i$.
The pairs $(k,\gamma)$ and $(l,\gamma')$ correspond to the same summands 
\[   f(i,k)\,\prod_{j=1\atop j\neq k}^{n+1} f(\gamma(j),j)  \, =\, f(i,l)\,\prod_{j=1\atop j\neq l}^{n+1} f(\gamma'(j),j) \]
-- yet, since $\sgn(\gamma)=-\sgn(\gamma')$, 
they have different sign and  cancel out each other.
\end{proof}


\end{document}